\documentstyle[11pt]{article}

\newcommand{\cD}{{\cal D}}

\newcommand{\cH}{{\cal H}}

\newcommand{\cK}{{\cal K}}

\newcommand{\B}[1]{\begin{#1}}
\newcommand{\E}[1]{\end{#1}}

\renewcommand{\P}{\mbox{$\bf P$}}     
\newcommand{\C}{\mbox{$\bf C$}}     
\newcommand{\Z}{\mbox{$\bf Z$}}     
\newcommand{\Q}{\mbox{$\bf Q$}}     

\newcommand{\G}{\mbox{$\Gamma$}} 
\newcommand{\D}{\mbox{$\Delta$}}

\newcommand{\df}{\mbox{\,$\stackrel{\pp{\rm def}}{=}$}\,}
\newcommand{\rat}{\mbox{$\stackrel{\pp{-\,-\succ}}{}$}}

\newcommand{\by}[1]{\stackrel{#1}{\rightarrow}}
\newcommand{\longby}[1]{\stackrel{#1}{\longrightarrow}}

\newcommand{\tensor}{\otimes}
\newcommand{\into}{\hookrightarrow}
 
\renewcommand{\iff}{\mbox{ $\Leftrightarrow$ }}

\newcommand{\ie}{{\it i.e.\/}\ }
\newcommand{\eg}{{\it e.g.\/}\ }
\newcommand{\cf}{{\it cf.\/}\ }

\newcommand{\p}[1]{\mbox{$\scriptstyle {#1}$}}   
\newcommand{\pp}[1]{\mbox{$\scriptscriptstyle {#1}$}}

\newcommand{\limdir}[1]{{\displaystyle{\mathop{\rm
lim}_{\buildrel\longrightarrow\over{#1}}}}\,}

\title{Balanced varieties\thanks{This Note is the
germ of a work in progress, I wrote it down for seminal
use; the first lecture I gave on the subject was on
February 16, 1996 at the Universitat de Barcelona,
Departament d'Algebra i Geometria, which I thank for the
invitation.}}

\author{by Luca {\sc Barbieri-Viale}\thanks{Member  of
GNSAGA of CNR, partially supported by ECC Science Plans}}

\date{}

\begin{document}

\maketitle

\B{abstract}
After the work of Bloch and Srinivas on correspondences
and algebraic cycles we begin the study of a birational
class of algebraic varieties determined by the property
that a multiple of the diagonal is rationally equivalent 
to a cycle supported on proper subschemes.
\E{abstract}

\section*{Introduction}

One of the main results due to Mumford \cite{MU} is the
observation that a complex projective non singular surface
having a non zero global holomorphic $2$-form has a non
zero Albanese kernel; the result was obtained after
Severi's work on algebraic cycles and holomorphic forms.
Because of Mumford's result the functor $A_0$ given by
$0$-cycles of degree zero modulo rational equivalence is
not representable, in contrast with the codimension $1$
case treated by Grothendieck and Murre. These facts lead to
the study of those varieties for which the Albanese kernel
is zero and to ``weak'' representability.\\
Bloch's proof \cite{Bl1} of the mentioned result by
Mumford reveals the motivic (read algebraic) nature of
the problem: he observed that the vanishing of the Albanese
kernel of a surface $X$ at a sufficiently large base field
extension yields the existence of 1-dimensional subschemes
$Z_1$ and $Z_2$, $2$-dimensional cycles $\G_1$
and $\G_2$ on $X\times X$, supported on $Z_1\times X$ and
$X\times Z_2$ respectively, such that some non-zero
multiple of the diagonal is rationally equivalent to
$\G_1 +\G_2$ over the base field. Then by investigating
the action of the correspondences $\G_1$, $\G_2$ on the
trascendental (read not algebraic) part of $\ell$-adic
\'etale cohomology one immediately obtains its vanishing.\\
The subsequent paper by Bloch and Srinivas \cite{BlS}
generalizes Bloch's argument to varieties of dimension $>
2$ showing its influence on algebraic, Hodge and arithmetic
cycles.\\ After all of that one has the temptation of
understanding the geometry of those varieties $X$ of pure
dimension $n$, for which we assume given proper closed
subschemes $Z_1$ and $Z_2$, $n$-dimensional cycles $\G_1$
and $\G_2$ on $X\times X$, supported on $Z_1\times X$ and
$X\times Z_2$ respectively, and a positive integer $N$
such that the following equation
$$N\D_X = \G_1 +\G_2$$
holds in the Chow group $CH_n(X\times X)$. I like to
call such a variety {\it balanced}\, because it appears
rather natural to regard $X$ ``balanced'' by $\G_1$ and
$\G_2$ on $Z_1$ and $Z_2$ in a motivic way as explained
above or as will become more clear in the following.\\ The
purpose of this Note is to begin the study (in Section~1)
of ``balanced varieties and balanced morphisms'' by showing
some basic geometric properties \eg balancing is a
birational property, stable under products; I also
explain the methods by Bloch and Srinivas \cite{BlS} (in
Section~2) obtaining applications in the ``language''
of balanced varieties; I'll sketch some examples (in
Section~3) and, finally by the way, I would draw a
picture of some expected properties of balanced
varieties.\\ I would like to express my best thanks to
P.Francia for his warm encouragement, his criticism and
his suggestions on some highlights of matters treated
herein.

\section{Basic geometry of balancing}

A variety will be an equidimensional reduced separated
scheme of finite type over a fixed base field $k$. A non
singular variety will be a regular scheme which is a
variety.

\B{defi} A balanced variety is a variety $X$ of dimension
$n$ over $k$ such that there exist {\it (i)}\, proper
closed subschemes $Z_1$ and $Z_2$ in $X$ which are
of finite type over $k$, {\it (ii)}\, $n$-dimensional
cycles $\G_1$ and $\G_2$ on $X\times_k X$ which are
supported on $Z_1\times_k X$ and $X\times_k Z_2$
respectively, and {\it (iii)}\, a positive integer $N$,
such that the following equation \B{equation}\label{bileq}
N\D_X = \G_1 +\G_2 \E{equation} holds in the Chow group
$CH_n(X\times_k X)$, where $\D_X$ is the canonical cycle
associated with the diagonal imbedding.\\ In this case we
will shortly say that we have a balancing of a variety. We
then will say that $X$ is balanced by $\G_1$ and $\G_2$ on
$Z_1$ and $Z_2$. We will say that $Z_1$ and $Z_2$ are the
balance {\it pans}. The {\it weight}\, of $X$ w.r.t.
$Z_1$ and $Z_2$ is the positive integer $\mbox{min}
\{\mbox{dim}Z_1,\mbox{dim}Z_2\}$, which we denote by
$w(X)$. 
\E{defi}

The basic nice properties of balanced varieties are the
following.

\B{prop}\label{local} Let $X$ be a variety of dimension
$n$. Let $U$ be a Zariski open dense subset of $X$ \eg  $X$
integral and $U\subset X$ any non-empty Zariski open. Then
$X$ is balanced {\em if and only if}\, $U$ is balanced.
\E{prop}  \B{proof} Since $U$ contains all generic points
we have that dim $Z\df X-U < n$ whence the diagonal cycle
$\D_X$ restricts to the diagonal cycle $\D_U$: so we may
restrict the balancing as well. Conversely, if $U$ is
balanced we have $N\D_U = \G_1^U+\G_2^U$ where $\G_1^U$ is
a cycle supported on $(Z_1\cap U)\times U$ and $\G_2^U$ is
a cycle supported on $ U\times (Z_2\cap U)$ for some $Z_1$
and $Z_2$ closed subschemes of $X$. Now we can lift
$\G_1^U$ to a cycle $\G_1+\zeta_1$ on $X\times X$ where
$\G_1$ is supported on $Z_1\times X$ and $\zeta_1$ is
supported on $(Z_1\cap Z)\times Z$; in the same way
$\G_2^U$ yields a cycle $\G_2+\zeta_2$ with $\G_2$
supported on $X\times Z_2$ and $\zeta_2$ supported on
$Z\times (Z_2\cap Z)$. Since the restriction of $\G_1+\G_2
+\zeta_1+\zeta_2$ to $U\times U$ is  $\G_1^U+\G_2^U$ and
the restriction of $N\D_X$ is $N\D_U$ we have that
$\G_1+\G_2 +\zeta_1+\zeta_2-N\D_X$ is a cycle supported on
$Z\times Z$. We then obtain an equation
$$N\D_X=\G_1+\G_2+\zeta$$  for some cycle $\zeta$
supported on $Z\times Z$, which is a balancing of
$X$.\E{proof} \B{prop}\label{cover} Let $X$ be a balanced
variety. Let $f: X\to X'$ be a proper surjective morphism
to a variety $X'$ with dim $X$ = dim $X'$. Then $X'$ is
balanced. \E{prop}
\B{proof} Since $f$ is generically finite we have that
the fundamental cycle $[X]$ has proper push-forward
$f_*[X]$ = deg($f$)$[X']$. Let denote by $\delta_X$ and
$\delta_{X'}$ the diagonal imbeddings of $X$ and $X'$. We
clearly have that $\delta_{X'*}\p{\circ}f_*=(f\times
f)_*\p{\circ}\delta_{X*}$ whence we have
$$(f\times f)_*(\D_X) = \mbox{deg}(f) \D_{X'}$$
by the definition of $\D$. Applying $(f\times f)_*$
to (\ref{bileq}) we have
$$N(f\times f)_*(\D_X) = (f\times f)_*(\G_1 + \G_2)$$
so that we obtain 
$$N\mbox{deg}(f) \D_{X'} = (f\times f)_*(\G_1) + (f\times
f)_*(\G_2)$$
Then $X'$ is balanced by $(f\times f)_*(\G_1)$ and 
$(f\times f)_*(\G_2)$ on $f(Z_1)$ and $f(Z_2)$.
\E{proof}

We then have:
\B{cor}\label{bir} Let $f:X'\rat X$ be a proper dominant
rational map with dim $X$ = dim $X'$: if $X'$ is balanced
then $X$ is balanced. Balancing is a birational property of
arbitrary varieties over arbitrary fields. \E{cor}
\B{proof} We may indeed restrict $f$ to open dense
subsets $U'$ of $X'$ and $U$ of $X$ such that $U'\to U$ is
a proper morphism and make use of the
Propositions~\ref{local},\ref{cover}. If $X$ is
birational to $X'$ then they have isomorphic open dense
subsets and we just use Proposition~\ref{local}.\E{proof}
\B{prop}\label{prod} Let $X$ be a balanced variety. Then
the product of $X$ with any variety is balanced.\E{prop}
\B{proof} Let $X'$ be any variety over $k$. 
We need the following simple Lemma.
\B{lemma} Let $X$ and $X'$ be varieties over $k$. Let
$$\sigma : X\times X \times X'\times X'\to X\times X' 
\times X\times X'$$ be the isomorphism given by 
$\sigma (x_1,x_2,x'_1,x'_2) = (x_1,x'_1,x_2,x'_2)$.
Then 
\B{equation}\label{simde}
\sigma_*(\D_X\times \D_{X'}) = \D_{X\times X'}
\E{equation}
\E{lemma}
\B{proof} Because of the functoriality of the exterior
product of cycles we have 
\B{equation}\label{extde}
(\delta_X\times\delta_{X'})_*\p{\circ}(-\times \dag) =
\delta_{X*}(-)\times\delta_{X'*}(\dag)
\E{equation}
Moreover we clearly have 
\B{equation}\label{compde}
\sigma\p{\circ}(\delta_X\times\delta_{X'}) =
\delta_{X\times X'}
\E{equation}
So then 
\B{center}
\parbox{7cm}{$\sigma_*(\D_X\times \D_{X'})=$\\
$=\sigma_*(\delta_{X*}[X]\times\delta_{X'*}[X'])=$\hfill
by (\ref{extde})\\
$=\sigma_*((\delta_X\times\delta_{X'})_*([X]\times
[X']))=$\\$=\sigma_*((\delta_X\times\delta_{X'})_*[X\times
X'])=$  \hfill by (\ref{compde})\\ $=(\delta_{X\times
X'})_*[X\times X']=$\\ $=\D_{X\times X'}$\\} 
\E{center}
as claimed.
\E{proof}

Let assume $X$ balanced by $\G_1$ and $\G_2$ on $Z_1$ and
$Z_2$ and let (\ref{bileq}) holds. We then have:
\B{center}
\parbox{7cm}{$N\D_{X\times X'}=$\hfill by (\ref{simde})\\
$=N\sigma_*(\D_X\times \D_{X'})=$\hfill
by (\ref{bileq})\\$=\sigma_*((\G_1
+\G_2)\times\D_{X'})=$\\$=\sigma_*(\G_1\times\D_{X'})
+\sigma_*(\G_2\times\D_{X'})$\\} 
\E{center}
where the support of $\sigma_*(\G_1\times\D_{X'})$ is
contained in $Z_1\times X'\times X\times X'$ and the
support of $\sigma_*(\G_2\times\D_{X'})$ is in 
$X\times X'\times Z_2\times X$. Then the
product $X\times_k X'$ is balanced on
$Z_1\times X'$ and $Z_2\times X'$. By symmetry also 
the product $X'\times_k X$ is balanced.
\E{proof}

\vspace{0.5cm}
We are now going to give another description of the
variation of a balancing under birational maps between
non-singular varieties. We recall that for any regular
algebraic scheme $V$ of pure dimension $n$ and $S\subset
V$ a closed subscheme we have a duality isomorphism
\B{equation}\label{kdual} H^p_S(V,\cK_p)\by{\simeq}
CH_{n-p}(S) \E{equation}
where $\cK_p$ is the Zariski sheaf on $V$ associated with
Quillen $K$-theory (\cf \cite{BV2}). If we have a
cartesian square
\B{displaymath}\B{array}{ccc}  
S' & \hookrightarrow & V' \\ 
 \downarrow & & \downarrow\p{f} \\
S & \hookrightarrow & V
\E{array}\E{displaymath}
we then have a commutative diagram
\B{displaymath}\B{array}{ccc}  
H^p_S(V,\cK_p) & \to & H^p(V,\cK_p) \\ 
\downarrow & & \downarrow \ \\
H^p_{S'}(V',\cK_p) & \to &  H^p(V',\cK_p)
\E{array}\E{displaymath}
whence by (\ref{kdual}) the following commutative diagram
of Chow groups
\B{equation}\label{ciao}\B{array}{ccc}  
CH_{n-p}(S) & \by{i_*} & CH_{n-p}(V) \\ 
\downarrow & & \downarrow\p{f^*} \\
CH_{n'-p}(S') & \by{{i'}_*} & CH_{n'-p}(V')
\E{array}\E{equation}
where $i$ and $i'$ are the imbedding of $S$ in $V$ and 
$S'$ in $V'$ respectively.\\ Let $f: X'\to X$ be a proper
birational morphism between non-singular varieties. For
$f: X'\to X$ as above we let consider the cartesian square
\B{displaymath}\B{array}{ccc}  
E_f& \by{\pi} & X \\ 
\p{e}\downarrow & & \downarrow\p{\delta} \\
X'\times X'& \by{f\times f} & X\times X
\E{array}\E{displaymath}
where $\delta$ is the diagonal imbedding. Because of
(\ref{ciao}) we obtain ($n =$ dim $X =$ dim $X'$)
\B{equation}\label{diciao}\B{array}{ccc}  
CH_{n}(X\times X) & \by{(f\times f)^*} & CH_{n}(X'\times
X') \\  \p{\delta_*}\uparrow & & \uparrow\p{e_*} \\
CH_{n}(X) & \by{{\pi}^*} & CH_{n}(E_f)
\E{array}\E{equation}
and moreover:\\[0.5cm]
{\it (i)}\, there is a closed imbedding
$i:X'\hookrightarrow E_f$ such that $e\p{\circ}i
=\delta_{X'}$ and $\pi\p{\circ}i =f$;\\[0.5cm]
{\it (ii)}\, if $Z'=f^{-1}(Z)$ is a closed in $X'$ such
that $f:X'-Z'\by{\simeq}X-Z$ then $e(E_f-X')\subset
Z'\times Z' - \D_{Z'}$;\\[0.5cm]
{\it (iii)}\, $E_f$ can also be regarded as the fiber
product of $f$ with itself \ie 
\B{displaymath}\B{array}{ccc}  
E_f& \by{\pi '} & X'\\ 
\p{\pi '}\downarrow & & \downarrow\p{f} \\
X'& \by{f} & X
\E{array}\E{displaymath}
is a cartesian square, where $f\p{\circ}\pi '=\pi$ and 
$\pi ':E_f\to X'$ has a section given by
$i:X'\hookrightarrow E_f$; we then have a splitting exact
sequence 
\B{equation}\label{split}
0\to CH_n(X')\by{\leftarrow}CH_n(E_f)\to CH_n(E_f-X')\to 0
\E{equation} because of ${\pi}'_*\p{\circ}i_* =
1$.\B{lemma} Notations and hypothesis as above. Then there
is a cycle $\zeta_f$ on  $X'\times X'$ supported on
$Z'\times Z'$ such that the following equation  
\B{equation}\label{zeta}
(f\times f)^*\D_{X}=\D_{X'}+\zeta_f
\E{equation}
holds in $CH_n(X'\times X')$.\E{lemma}
\B{proof} Because of (\ref{diciao}) we have that 
$$(f\times f)^*\D_{X}=(f\times f)^*(\delta_{X*}[X])
=e_*\pi^*[X]$$
So we are left to show the equation
$e_*\pi^*[X]=\D_{X'}+\zeta_f$. The projection ${\pi}'_*$ of
$\pi^*[X]\in CH_{n}(E_f)$ is $[X']$ since the composite of
$$H^n_X(X\times X,\cK_n)\to H^n_{E_f}(X'\times X',\cK_n)
\to H^n_X(X\times X,\cK_n)$$ is the identity, by
\cite[Section~7]{BV2} because $f\times f$ is
birational, and the composite of
$$CH_{n}(E_f)\by{{\pi}'_*}CH_n(X')\by{f_*}CH_n(X)$$ is
$\pi_*$ by {\it (iii)}\, above, where $f_*$ is the
isomorphism sending $[X']$ to $[X]$ because $f$ is mapping
birationally each component of $X$ to a distinct component
of $X'$.\\ Moreover by {\it (ii)}\, we have a commutative
diagram $$\begin{array}{cc}
CH_{n}(E_f-X') \ \ \to \ \ CH_{n}(Z'\times Z'
-\D_{Z'}) & {\by{\simeq}}\ \ CH_{n}(Z'\times Z') \\
 \p{e_*}\searrow \hspace*{40pt} {\downarrow} 
&{\downarrow}\\ \ \hspace*{80pt} CH_{n}(X'\times
X'-\D_{X'})  & {\leftarrow}\ \ \  CH_{n}(X'\times X')  
\end{array}$$
since $Z'$ has dimension $\leq n-1$ whence $CH_{n}(Z')=0$.
Thus we are able to give a description of $e_*\pi^*[X]$
w.r.t.\/ the splitting (\ref{split}): from the latter
diagram we see that the restriction of $\pi^*[X]$ to
$CH_{n}(E_f-X')$ (\ie a component of $\pi^*[X]$) has
image (\ie $e_*$ of it in  $CH_{n}(X'\times X')$ is) a
cycle $\zeta_f$ supported on $Z'\times Z'$; the component
of $\pi^*[X]$ in $CH_{n}(X')$ (\ie $i_*[X']$) has image
$\D_{X'}$ because of {\it (i)}\, \ie $e_*\p{\circ}i_*
=\delta_{X'*}$. The claimed formula (\ref{zeta}) is
obtained.\E{proof}

Assuming $X$ balanced we have:
\B{center}
\parbox{7cm}{$(f\times f)^*(\G_1 +\G_2)=$\hfill
by (\ref{bileq})\\   
$=N(f\times f)^*\D_{X}=$\hfill
by (\ref{zeta})\\$=N\D_{X'}+N\zeta_f$\\} 
\E{center}
yielding the equation 
$$N\D_{X'} = (f\times f)^*(\G_1)+(f\times
f)^*(\G_2) -N\zeta_f$$
Now, because of (\ref{kdual}) and (\ref{ciao}), we have
that the cycle  $(f\times f)^*(\G_1)$ is supported on 
$f^{-1}(Z_1)\times X'$ and $(f\times f)^*(\G_2) -N\zeta_f$
on $X'\times f^{-1}(Z_2)\cup Z'\times Z'$. Then $X'$ is
balanced on $f^{-1}(Z_1)$ and
$f^{-1}(Z_2)\cup Z'$.\\[0.5cm]
{}\hfill
In the following any flat morphism will have
a relative dimension.

\B{defi} Let $p:X\to S$ be a morphism of
varieties over $k$ where $S$ is irreducible of dimension
$i$. We say that $p:X\to S$ is {\it balanced}\, or that $X$
is balanced {\it over}\, $S$ if $p:X\to S$ is
flat of relative dimension $n$ and the equation
$\Delta_{X/S}= \G_1+\G_2$ holds in $CH_{n+i}(X\times_SX)$
where $\G_1$ is supported on  $Z_1\times_SX$ and $\G_2$ is
supported on $X\times_SZ_2$ for some $Z_1$ and $Z_2$
proper closed subschemes of $X$ which are flat over $S$.
\E{defi} \B{prop} We have the following basic properties.
\B{description} \item[{\sc Base change}] Let $p:X\to S$ be
a balanced morphism, $\varphi :S'\to S$ be any flat
morphism (but $S'$ irreducible) and $X'= X\times_SS'$.
Then $p':X'\to S'$, obtained by base extension, is
balanced.\\ \item[{\sc Composition}] If $p:X\to S$ is
balanced and $\varphi :S\to T$ is flat ($T$ irreducible)
then $\varphi\p{\circ}p :X\to T$ is balanced.\\ \item[{\sc
Product}] If $X$ is balanced over $S$ and $Y$ is flat over
$S$ then $X\times_SY$ is balanced over $S$.
\E{description}
\E{prop} \B{proof} A balancing for the base
extension $p':X'\to S'$ can be obtained by making use of
external products of flat families of cycles because of
the equation $\D_{X'/S'}=\D_{X/S}\times_{S}[S']$ (\cf the
proof of Proposition~\ref{prod}). The composition
$\varphi\p{\circ}p :X\to T$ is known to be flat and the
imbedding $j:X\times_SX\into X\times_TX$ is known to be
closed, yielding the equation $j_*\D_{X/S}=\D_{X/T}$.
Finally: $X\times_SY\to S$ is the composition of the flat
morphism $Y\to S$ with the base extension.
\E{proof}

By the above, if $X$ is balanced over $S$ then $X$ is
balanced over the base field. Moreover $X$ can be view as
an algebraic family and, for example, we have that $X_s$, 
the fibre at the generic point $s\in S$, is balanced over
$K(S)$ (the function field of $S$). Indeed, let
$V$ be any Zariski open neighborhood of the generic point
$s$ in $S=\overline{\{s\}}$: by flat base change we have
that $X\times_SV\to V$ is balanced. The claim above is
obtained by taking the limit over $V$ because we have
$$CH_n(X_s)\cong \limdir{V} CH_n(X\times_S V)$$ and
$\D_{X_s/K(S)} =\limdir{V}\D_{X\times_S V/V}$.\\

\B{defi}\label{universal} Let $X$ be a variety over a field
$k$. We say that $X$ is $K$-balanced if $X_K$ is balanced
over the field extension $k\subset K$. If $\Omega$ is a
universal domain in the sense of Weil and $X$ is
$\Omega$-balanced we say that $X$ is {\it universally}\,
balanced. \E{defi} We will see in the next Section, that
universally balanced varieties are also balanced over the
base field. 

\section{Actions of a balancing}

Let $(H^*,H_*)$ be an appropriate duality theory in
the sense of \cite[Sections 6-7]{BV2} \eg De Rham theory
or  $\ell$-adic \'etale theory. The
twisted cohomology functor $H^*(-,\cdot)$ is then
equipped with a functorial `cycle class' map
$c\ell:H^p(X,\cH^p(\p{p})\to H^{2p}(X,\p{p})$ for $X$
smooth over $k$. An algebraic `$\cH$-correspondence' from
$X$ to $Y$ is an $\cH$-cohomology class $\alpha\in
H^{n+r}(X\times Y, \cH^{n+r}(\p{n+r}))$ ($X$ has pure
dimension $n$) and it will be denoted by $\alpha :
X\leadsto Y$. As usual $\alpha$  acts on the cohomology
groups of non-singular projective
varieties $$\alpha_{\sharp}: H^i(X,\p{j})\to
H^{i+2r}(Y,\p{j+r})$$ where $\alpha_{\sharp}(-)\df
p_{Y,*}(c\ell(\alpha)\p{\cup} p^*_X(-))$ and $p_{Y}$,
$p_{X}$ are the projections of $X\times Y$ on $Y$ and $X$
respectively. Because of our assumptions there is a
canonical ring homomorphism from the Chow ring to the
$\cH$-cohomology ring (\cf \cite[Sections 5.5 and
6.3]{BV2}). The action of any algebraic cycle on
a fixed cohomology theory is given {\it via}\, its 
$\cH$-cohomology incarnation \eg any cycle 
algebraically equivalent to zero acts as zero on the
singular cohomology of the associated analytic space. In
particular any  $\alpha : X\leadsto Y$ as above acts on the
$\cH$-cohomology groups $$\alpha_{\sharp}:
H^p(X,\cH^q(\p{j}))\to H^{p+r}(Y,\cH^{q+r}(\p{j+r}))$$ in
a compatible way w.r.t. the coniveau spectral sequences.\\
An easy application of the projection formula (for
$\cH$-cohomologies we need \cite[Section 5.4]{BV2}) yields
the following useful Lemma. \B{lemma} Let $f: X'\to X$ and
$g: Y'\to Y$ be morphisms of non-singular projective
varieties. For $\alpha: X'\leadsto Y$ we have
\B{equation}\label{uno}
(f\times 1_Y)_*(\alpha)_{\sharp} = \alpha_{\sharp}
\p{\circ}f^* \E{equation}
For $\beta : X\leadsto Y'$ we have  
\B{equation}\label{due}(1_X\times g)_*(\beta)_{\sharp}=
g_*\p{\circ}\beta_{\sharp}
\E{equation}
\E{lemma}
\B{proof} The proof is left as an exercise for the reader.
\E{proof}

Let now consider a balanced variety $X$ which is
non-singular and projective; we are going to make use of
(\ref{bileq}) w.r.t. the various cohomology theories. We
need to assume resolution of singularities and we let
assume that $Z_1$ and $Z_2$ are equidimensional. We then
may consider a resolution $Z_1'\to Z_1$ and we have that
$\G_1 = (f\times 1_X)_*(\G_1')$ where $\G_1'$ is a cycle on
$Z_1'\times X$ of codimension equal to the dimension of
$Z_1$ and $f : Z_1'\to X$, so that by (\ref{uno}) we 
obtain  \B{equation} (\G_1)_{\sharp} = (f\times
1_X)_*(\G_1')_{\sharp}= (\G_1')_{\sharp}\p{\circ}f^*
\E{equation}
since $\G_1': Z_1'\leadsto X$; \ie we have a commutative
triangle
\B{eqnarray}\label{map1}
H^p(X,\cH^q\p{(\cdot)}) & \longby{\G_{1\sharp}} &
H^p(X,\cH^q\p{(\cdot)}) \nonumber \\
 \searrow    &  & \nearrow \\
     &  H^p(Z_1',\cH^q\p{(\cdot)}) &\nonumber
\E{eqnarray}
We consider as well a resolution of singularities
$Z_2'\to Z_2$ and we have $\G_2 = (1_X\times g)_*(\G_2')$
where $\G_2'$ is a cycle on $X\times Z_2'$ whence the
degree of the correspondence $\G_2':X\leadsto Z_2'$ is
$-\mbox{codim}_X(Z_2)$; by (\ref{due}) we then have
\B{equation}
(\G_2)_{\sharp} = (1_X\times g)_*(\G_2')_{\sharp} =
g_*\p{\circ}(\G_2')_{\sharp}
\E{equation}
\ie the following triangle
\B{eqnarray}\label{map2}
H^p(X,\cH^q\p{(\cdot)}) & \longby{\G_{2\sharp}} &
H^p(X,\cH^q\p{(\cdot)}) \nonumber \\
 \searrow    &  & \nearrow \\
     &  H^{p-c}(Z_2',\cH^{q-c}\p{(\cdot-c)}) &\nonumber
\E{eqnarray}
commutes, where $c=\mbox{codim}_X(Z_2)$. If $Z_2$ (or
$Z_1$) is not equidimensional we may consider each smooth
component of its resolution $Z_2'$ acting as above.
\B{prop} Let consider any cohomology theory $H^*$ as
above. Let cd$(k)$ be the ``cohomological dimension'' of
the field $k$ \ie we assume that $H^i(U) = 0$ if
$i>\mbox{dim} U +\mbox{cd}(k)$ and $U$ is affine. If $X$
is balanced of weight $w$, and either the pans are
smooth or can be resolved, then
$H^0(X,\cH^q\p{(\cdot)})$ is {\em $N$-torsion} for
$q>w+\mbox{cd}(k)$ and some positive integer $N$.
\E{prop}  
\B{proof} Is essentially the same of \cite[Th.
1, ii--iii]{BlS}. By interchanging the pans we may assume
that $w=$dim$Z_1$. The action of the cycle $N\D_X$ is the
multiplication by $N$ which is given  by
$\G_{1\sharp}+\G_{2\sharp}$. Beacause of the assumptions
and (\ref{map1}) we have $\G_{1\sharp}=0$   since
$\cH^q=0$ on $Z_1'$ if $q>w+\mbox{cd}(k)$. Because  of
(\ref{map2}) $\G_{2\sharp}=0$ since
$\mbox{codim}_X(Z_2)>0$. 
\E{proof}

Whenever the sheaves $\cH^*\p{(\cdot)}$ are torsion free
we have that $H^0(X,\cH^*\p{(\cdot)})=0$ under the
assumptions in the Proposition above. Because of \cite{BV1}
(\cf \cite{BV0}, \cite{BV3}) we then have: \B{cor} Let $X$
be balanced over $\C$ of weight $w$ and let assume that
Kato's conjecture hold true for function fields. If $q>w$
we then have $$H^0(X,\cH^q(\Z))=0$$ and
$$H^0(X,\cH^{q+1}(\Z(q)_{\cD}))=0$$
where $\cH^*(\Z)$ and $\cH^{*}(\Z(\cdot)_{\cD})$ are the
Zariski sheaves associated with singular cohomology and 
Deligne--Beilinson cohomology respectively. 
\E{cor}
As in \cite[Th. 1]{BlS} by the above one obtains the
following: for $X$ balanced over $\C$ of weight $\leq 2$ 
algebraic and homological equivalence coincide for
codimension $2$ cycles and for $X$ of weight $\leq 1$
also homological and Abel-Jacobi equivalence coincide in
codimension $2$. Indeed the Griffiths group is a
quotient of $H^0(X,\cH^3(\Z))$ and the Abel-Jacobi kernel
is a quotient of the torsion
free group $H^0(X,\cH^{3}(\Z(2)_{\cD}))$ (\cf
\cite{BV3} for surfaces).\\   
Moreover:
\B{prop}\label{forms} Let $X$ (projective, non-singular) be
balanced over $\C$ of weight $w$. If $q>w$ then
$$H^0(X,\Omega^q_X)=0$$ 
\E{prop} 
\B{proof} The previous arguments applies to the twisted
cohomology theory given by $F^{\cdot}H^*$ \ie the De Rham
filtration, and $F^qH^q(X)=H^0(X,\Omega^q_X)$.  Now
$\G_{1\sharp}=0$ because $\Omega^q_{Z_1'}=0$ if $q>w=$ dim
$Z_1$ and (\ref{map1}). Moreover $\G_{2\sharp}=0$ by
$F^{q-c}H^{q-2c}(Z_2')=0$ since $c>0$ \ie the codimension
of any component of $Z_2$ cannot be zero. Since 
$F^{\cdot}H^*$ are real vector spaces we are done.
\E{proof}
\B{prop}\label{Alba} If $X$ (projective, non-singular)
is balanced of weight $w$ over a field $k$, and either the
pans are smooth or can be resolved, then there exists a
closed  subscheme $Y$ of $X$ such that $CH^i(X-Y)_{\Q}=0$
for all $i>w$. For $0$-cycles $Y$ exists of dimension
$w$. If, moreover, $k$ is algebraically closed then
$A^i(X-Y)=0$ for all $i>w$ and $A_0(X-Y)=0$ as above;
finally, in this case, the Albanese kernel of $X$ is
contained in the Albanese kernel of the resolved pan of
dimension $w$. \E{prop} \B{proof} We may assume that $Z_1$
is irreducible and non-singular of dimension $w$:
therefore we obtain that $\G_{1\sharp}=0$ on $CH^i(X)$ if
$i>w$, by (\ref{map1}) and $CH^i(Z_1)=0$. Then 
$\G_{2\sharp}\tensor\Q$ is the isomorphism induced by the
multiplication by $N$ on $CH^i(X)_{\Q}$ if $i>w$; then, by
(\ref{map2}), the Chow group of $Z_2$ surjects onto 
$CH^i(X)_{\Q}$ if $i>w$. By taking $Z_2$ as $Y$ we have
that $CH^i(X-Y)_{\Q}=0$ as claimed. If $i=\mbox{dim} X$
\ie for $0$-cycles, we then can take $Z_2$ of dimension
$w$ and we conclude as above. If $k$ is closed it is well
known that $A_0(X)$ is divisible and the Albanese kernel
is uniquely divisible; by the same argument we obtain the
claimed results. \E{proof}

Following Bloch, Srinivas and Jannsen (\cf
\cite{Bl1}, \cite{BlS}, \cite[Remark~1.7]{JA}) we
have:
\B{prop} \label{rap} Let $X$ be a non-singular projective
variety over a field $k$. Then $X$ is universally balanced
of weight $\leq 1$ if and only if the Chow group of
$0$-cycles of degree zero is representable. \E{prop}
\B{proof} By the definition, $X$ is
universally balanced if $X_{\Omega}$ is balanced over a
universal domain $\Omega$. By the Proposition~\ref{Alba}
the Albanese kernel $X_{\Omega}$ is contained in the
Albanese kernel of (at most) a smooth curve over $\Omega$,
which is zero. Conversely, if the Chow group of $0$-cycles
is representable then by \cite[Prop.~1.6]{JA} there is
(at most) a smooth projective curve over $\Omega$  and a
morphism $g:C\to X_{\Omega}$ such that
$A_0(X_{\Omega}-g(C))=0$ hence by \cite[Prop.~1]{BlS}
$X_{\Omega}$ is balanced of weight $\leq 1$. \E{proof}
\B{prop} Let $X$ be a non-singular projective variety over
a field $k$. If $X$ is universally balanced then it is
balanced. \E{prop} \B{proof}
The proof is similar to \cite[Proof of
Th.~3.5.(a)-(b)]{JA} which is quite the same of the proof
of Proposition~1 of \cite{BlS}.\E{proof}

Following Jannsen \cite[\S~3]{JA} we have:  
\B{cor} If $X$ is universally balanced of weight $\leq w$
then $H^i(X,\p{j})$ (= any twisted cohomology theory in the
sense of \cite[\S~3]{JA}) is of coniveau $1$ for
$i>w$.\E{cor} Finally, by Proposition~\ref{forms} and
Proposition~\ref{rap} we can easily obtain the following.
\B{cor} {\em (Mumford-Roitman Theorem)} If
$H^0(X,\Omega^q_X)\neq 0$ for some $q>1$ then the Chow
group of $0$-cycles of degree zero is not representable.
\E{cor}
See also \cite[Cor.~3.7]{JA}.

\section{Paradigma} We now give some examples.\\

\B{example} The projective space $\P^n_k$ over any field
$k$ is balanced of weight $0$ \ie the pans are given by
any couple of closed points $x_0$ and $x_1$. In fact  we
have that 
$$CH_n(\P^n_k\times \P^n_k) \cong (x_0\times
\P^n_k)\Z\bigoplus (\P^n_k\times x_1)\Z$$
\E{example}
\B{example} Because of Corollary~\ref{bir} and the above,
unirational varieties are balanced over any field. More
generally, because of Proposition~\ref{prod}, we have
that uniruled varieties are balanced, since are dominated
by a product with $\P^1_k$.
\E{example}
\B{example} If $X$ is a smooth closed subvariety of an
abelian variety then $X$ is not balanced over $\C$, since
$H^0(X,\Omega^n_X)\neq 0$ for $n=$ dim $X$ and
Proposition~\ref{forms}.
\E{example}
\B{example} The Kummer $3$-fold over $\C$ is balanced of
weight $2$ by Bloch and Srinivas \cite{BlS} since $A_0(X)$
is generated by cycles supported on a finite number of
surfaces. This is an example of balanced variety which is
quite far from ruled varieties.
\E{example}

\section*{Some remarks and questions}

\B{enumerate}

\item Of course one might like to define balanced
schemes starting from an equidimensional separated scheme
over an arbitrary base scheme but then does one have
anything to say about it ? Nevertheless one case that
sounds interesting is the case of arithmetic schemes \ie
regular schemes which are projective and flat over the
integers, by mean of arithmetic cycles and correspondences
in the sense of H. Gillet and C. Soul\'e. So I expect a
parallel theory of ``arithmetically balanced varieties'' by
analysing the action of arithmetic correspondences.

\item In the language of balanced varieties, Bloch's
conjecture give us a numerical criterion for balanced
surfaces \ie $p_g=0$ for weight $1$ and $p_g=q=0$ for
weight zero. We may expect that the following numerical
criterion $$X\mbox{ balanced of weight }w \iff
H^0(X,\Omega^q_X) = 0 \mbox{ for } q>w$$
holds true for any projective complex manifold $X$.\\
More in general, for $X$ defined over any field $k$,
$X$ will possibly be universally balanced of weight $\leq
w$ if and only if any suitable `cohomology theory' will be
of coniveau $1$ in degrees $>w$.

\item Since the global forms are invariants under a smooth
deformation, the previous speculation is suggesting us
that the balancing should be a deformation property. We
remark that uniruled varieties do have this property,
thanks to A. Fujiki and M. Levine (deformation invariance
of rationally connected varieties is due to J. Kollar, Y.
Miyaoka and S. Mori).

\item Let say that a flat family of varieties, 
parametrized by a nice variety, is balanced if
its general member is a balanced variety. Is there a
section ? 

\E{enumerate}

\vspace{1cm}

\B{flushright}
Dipartimento di Matematica\\
Universit\`a di Genova\\
Via Dodecaneso, 35\\
16146 -- {\sc Genova}\\
{\it Italia}
\E{flushright}


\begin{thebibliography}{}

\bibitem{BV0}{\sc L.Barbieri-Viale}: Des invariants
birationnels associ\'es aux th\'eories
cohomologiques, {\it C.R. Acad. Sci. Paris} {\bf 315}
S\'erie I (1992), 1259-1262.

\bibitem{BV1}{\sc L.Barbieri-Viale}: On the
Deligne--Beilinson cohomology sheaves, unpublished, 1994.

\bibitem{BV2}{\sc L.Barbieri-Viale}: $\cH$-cohomologies
versus algebraic cycles, preprint, 1994, to appear in
Math. Nachr., 51 p.

\bibitem{BV3}{\sc L.Barbieri-Viale}~and{\sc~V.Srinivas}: 
A reformulation of Bloch's conjecture,
{\it C.R. Acad. Sci. Paris} {\bf 321} S\'erie I (1995),
211-214. 

\bibitem{Bl1}{\sc S.Bloch}: Lectures on algebraic
cycles, {\it Duke Univ. Math. Series} {\bf 4}, Durham,
1980. 

\bibitem{BlS}{\sc S.Bloch}~and{\sc~V.Srinivas}: Remarks on
correspondences and algebraic cycles, {\it Amer. J.
Math.} {\bf 105} (1983), 1235-1253. 

\bibitem{JA}{\sc U.Jannsen}: Motivic sheaves and
filtrations on Chow groups, in Motives {\it AMS Proc. of
Symp. in Pure Math.} Vol. {\bf 55} Part 1, 1994.

\bibitem{MU}{\sc D.Mumford}: Rational equivalence of
$0$-cycles on surfaces, {\it J.Math.Kyoto Univ.} {\bf 9}
(1968) 195--204.

\end{thebibliography}
\end{document}